\documentclass[12pt]{article}

\usepackage{amssymb,bbm}    
\usepackage{amsmath}    
\usepackage{amsthm} 
\usepackage{mathtools}
\usepackage{mathrsfs}
\usepackage[x11names]{xcolor}
\usepackage{hyperref}

\hypersetup{colorlinks=true,urlcolor=black,linkcolor=black,citecolor=black,naturalnames=true,hypertexnames=false}
\usepackage[backend=bibtex,giveninits, style=authoryear, uniquelist=false, natbib=true, doi=false, url=false, isbn = false,maxcitenames=4,maxbibnames=99]{biblatex}

\usepackage{bib}

\usepackage{macros}

\theoremstyle{plain}
\newtheorem{theorem}{Theorem}[section]

\theoremstyle{remark}
\newtheorem{remark}[theorem]{Remark}

\makeatletter
\def\@fnsymbol#1{\ensuremath{\ifcase#1\or \star \or \ddag\or \S\else\@ctrerr\fi}}
\makeatother

\title{On rate optimal private regression under local differential privacy}
\author{L\'aszl\'o Gy\"orfi\thanks{Department of Computer Science and Information Theory, Budapest University of Technology and Economics, Magyar Tud\'{o}sok krt. 2., Budapest, H-1117, Hungary.
\break
\texttt{gyorfi@cs.bme.hu}}
\and Martin Kroll\thanks{Ruhr-Universit\"at Bochum, Fakult\"at f\"ur Mathematik, Universit\"atsstra\ss e 150, D-44801 Bochum.} \thanks{Universität Bayreuth, Fakultät für Mathematik, Physik und Informatik, D-95440 Bayreuth.  \texttt{martin.kroll@uni-bayreuth.de}.}}

\begin{document}

\maketitle

\begin{abstract}
We consider the problem of estimating a regression function from anonymized data in the framework of local differential privacy.
We propose a novel partitioning estimate of the regression function, derive a rate of convergence for the excess prediction risk over Hölder classes, and prove a matching lower bound.
In contrast to the existing literature on the problem the so-called strong density assumption on the design distribution is obsolete.
\end{abstract}

\noindent


\noindent


\section{Introduction}

Let $(X,Y) \in \R^d \times \R$ be a pair of random variables with explanatory variable $X \in \R^d$ and real-valued response $Y$ satisfying $\EXP[Y^2]<\infty$.
We denote by $\mu$ the distribution of $X$, that is, $\mu(A)=\PROB (X\in A)$ for all measurable sets $A\subseteq \R^d$.
Then, the \emph{regression function}
\begin{equation*}
\label{def:model}
 m(x) = \EXP [Y | X=x]
\end{equation*}
is well-defined for $\mu$-almost all $x \in \R^d$.
For any measurable function ${g\colon\R^d\to\R}$ one has
\[
\EXP [ ( g(X) - Y)^{2}] = \EXP  [( m(X) - Y )^{2} ]  + \EXP  [( m(X) - g(X) )^{2} ],
\]
and therefore, setting
\[
L^*=\EXP  [( m(X) - Y )^{2} ],
\]
it follows that
\begin{equation*}
\EXP [ ( g(X) - Y)^{2}] = L^*  + \int ( m(x) - g(x) )^{2} \mu(\dd x).
\end{equation*}
Measuring the performance of an estimator $\mhat$ of $m$ through the loss function
\[
L(m, \mhat) \defeq \int (m(x) - \mhat(x))^2 \mu(\dd x)
\]
may thus be interpreted as the excess prediction risk at a new design point $X$ distributed according to the design measure $\mu$.

In this paper, we consider piecewise constant estimators of the regression function $m$ based on cubic partitions.
Let $\P_h=\{A_{h,1},A_{h,2},\ldots\}$ be such a cubic partition
of $\Rd$ with cubic cells $A_{h,j}$ of volume
$h^d$.
The raw data $\mathcal D_n$ are assumed to be i.i.d.\,copies of the random vector $(X,Y)$,
\begin{align}\label{eq:raw_data}
\mathcal{D}_n \defeq \{(X_1,Y_1),\ldots,(X_n,Y_n)\}.
\end{align}
Put
\begin{align*}
\nu_n(A_{h,j})=\frac 1n \sum_{i=1}^n Y_i \IND_{\{X_i \in A_{h,j}\}}
\end{align*}
and
\begin{align*}
\mu_n(A_{h,j})=\frac 1n \sum_{i=1}^n \IND_{\{X_i \in A_{h,j}\}}.
\end{align*}
Then, a standard regression estimate is defined by
\[
m_{n}(x)= \frac{ \nu_n(A_{h,j}) }{ \mu_n(A_{h,j})}
\]
for any $x \in A_{h,j}$ with the usual convention that $0/0 \defeq 0$.
Theorem~4.3 in the monograph \citep{GyKoKrWa02} states an upper bound on the rate of convergence for this partitioning estimate for Lip\-schitz continuous regression functions.
Its extension to the more general case of H\"older continuous functions is straightforward, and yields the following result.
\begin{theorem} 
\label{brate}
If the function
\begin{align*}
\sigma^2(x) \defeq \Var(Y|X=x)
\end{align*}
is bounded, $m$ is $(\beta,C)$-H\"older smooth with index $0<\beta\le 1$, that is,
\begin{align}
\label{eq:hoelder}
|m(x)-m(x^\prime)|\le C\| x-x^\prime\|^{\beta} , \quad x,x^\prime\in \R^d,
\end{align}
and $X$ is bounded, then
\begin{equation}\label{eq:upper:nonprivate}
\EXP \left[ \int ( m(x) - m_{n}(x) )^{2} \mu(\dd x) \right]
 \lesssim \frac{1}{n {h}^d} + {h}^{2\beta}.
\end{equation}
\end{theorem}
In particular, the choice $h = h_n \asymp n^{-1/(2\beta+d)}$ realizes the best compromise between the two antagonistic terms on the right-hand side of \eqref{eq:upper:nonprivate}, and the resulting rate is $n^{-2\beta/(2\beta+d)}$.
Standard arguments for nonparametric lower bounds, for instance based on Assouad's lemma, show that this rate is indeed optimal.

\medskip

The main purpose of this paper is to give an analogue of Theorem~\ref{brate} for the case when the raw data $\mathcal D_n$ are not directly accessible but only a suitably anonymized surrogate.
More precisely, the anonymized data must satisfy a \emph{local differential privacy} (LDP) condition.
Our work is motivated by the recent paper \citep{BeGyWa21} where a first step in this direction was done.
In that paper the authors considered a private partitioning estimate and derived the upper bound $n^{-1/(d+1)}$ on the rate of convergence for Lipschitz continuous functions ($\beta = 1$).
However, this rate was established under a quite restrictive assumption on the design distribution $\mu$ (called \emph{strong density assumption} (SDA) in \citep{BeGyWa21}).
Moreover, it was conjectured that the rate of convergence can be arbitrarily slow when the SDA is not fulfilled.
In this paper, we will show that this conjecture does not hold.
Quite the contrary, we will even introduce an estimator that attains the rate $n^{-\beta/(\beta+d)}$ without the SDA.
We emphasize that even the existence of a Lebesgue density for $\mu$ will not be postulated.
We complement our upper bound by proving a minimax lower bound showing that the rate $n^{-\beta/(\beta + d)}$ is indeed optimal which is in accordance with what can be expected from similar problems like nonparametric density estimation where a similar deterioration of the rate of convergence has been found \citep{duchi2018,Butucea2020}. The same phenomenon has also been observed for classification problems \citep{BerrettButucea2019}.

The rest of the paper is organized as follows: For the sake of self-contained\-ness we recap the notion of LDP in Section~\ref{s:privacy} and introduce a suitable anonymization of the raw data that generates locally differentially private data.
In Section~\ref{Sec:MainRegression} we introduce a modification of the classical partitioning estimate of the regression function that is only based on the availability of the anonymized data and derive a convergence rate for this estimator.
In Section~\ref{sec:lower} we prove a matching lower bound that coincides with the upper one.
All proofs are gathered in Section~\ref{Sec:Proofs}.

\section{Anonymization of the raw data}
\label{s:privacy}

In this section, we briefly recall the definition of LDP before we describe our privacy mechanism.
In the language of probability theory, \emph{non-interactive} privacy mechanisms are given by conditional distributions $Q_i$, $i=1,\ldots,n$, that draw privatized data $Z_i$ from potentially different measurable spaces $(\Zc_i,\Zs_i)$.
More precisely, given the raw data $(X_i,Y_i)=(x_i,y_i)$, one draws $Z_i$ according to a probability measure defined through $Q_i(A|(X_i,Y_i)=(x_i,y_i))$ for any $A \in \Zs_i$.
Such a non-interactive mechanism is \emph{local} since any data holder can independently generate privatized data.

For a privacy parameter $\alpha \in [0,\infty]$, any non-interactive privacy mechanism is said to be an $\alpha$-\emph{locally differentially private} mechanism if the condition
\begin{equation}\label{eq:LDP:noninteractive}
    \frac{Q_i(A|(X_i,Y_i)=(x,y))}{Q_i(A|(X_i,Y_i)=(x\p,y\p))} \leq \exp(\alpha)
\end{equation}
is satisfied for any $A \in \Zs_i$ and all potential values $(x,y),(x\p,y\p)$ of the raw data.
The set of all $\alpha$-locally differentially private mechanisms will be denoted by $\Qc_\alpha$.

\bigskip

Let us now state the specific privacy mechanism that we consider in this work for the anonymization of the raw data $\mathcal D_n$ in \eqref{eq:raw_data}.
Our approach follows the technique of Laplace perturbation already considered in \citep{duchi2018,BerrettButucea2019,Butucea2020,BeGyWa21,Gyoerfi2023}.
In order to define the privatized data, we first choose a closed Euclidean ball $B = \{ x \in \R^d : \lVert x \rVert \leq r \}$ of radius $r>0$ centred at the origin.
Let further $A_{h,1},A_{h,2},\ldots$ be a partition of $\R^d$ consisting of cubic cells with volume $h^d$ for some  $h=h_n > 0$.
Without loss of generalization we can assume that the cells are numbered such that $A_{h,j}\cap B \ne \emptyset$ when $j\le N_n$ for some non-negative integer $N_n$, and $A_{h,j}\cap B= \emptyset$ otherwise.
For a threshold $T > 0$ and any $x \in \R$ we define $[x]_T = \max\{ -T, \min \{ x, T \} \}$.
In our privacy setup, the data holder of the $i$-th datum $(X_i,Y_i)$ generates and transmits to the statistician the data
\begin{equation}
\label{Eq:Mech1}
	Z_{ij} \defeq  [Y_i]_{T} \IND_{\{X_i \in A_{h,j}\}} + \sigma_Z \epsilon_{ij}, \quad j = 1,\ldots,N_n,
\end{equation}
and
\begin{equation}
\label{Eq:Mech2}
W_{ij} \defeq \IND_{\{X_i\in A_{h,j} \}} + \sigma_W \zeta_{ij}, \quad j = 1,\ldots,N_n,
\end{equation}
with noise levels $\sigma_Z, \sigma_W>0$ and $\epsilon_{ij}$, $\zeta_{ij}$ ($i=1,\ldots,n$, $j=1,\ldots,N_n$) are independent centred Laplace random variables with unit variance.
This means that the individual with index $i$ generates noisy data for any cell $A_{h,j}$ that has non-trivial intersection with the ball $B$.
The noise levels $\sigma_Z$ and $\sigma_W$ have to be chosen sufficiently large in dependence on the desired privacy level $\alpha$ in order to make the overall mechanism satisfy $\alpha$-LDP.
It has been shown in \citep{BeGyWa21}, p.~2438, that the choices
\begin{equation}\label{eq:choice:sigmas}
    \sigma_W^2 = 32/\alpha^2 \quad \text{and} \quad \sigma_Z^2 = 32T^2/\alpha^2
\end{equation}
ensure $\alpha$-LDP.

\begin{remark}\label{rem:randomized_response}
The privacy mechanism defined by \eqref{Eq:Mech1} and \eqref{Eq:Mech2} is certainly not the only potential privacy mechanism.
For instance, one could build on the randomized response technique discussed in \citep{duchi2018}.
We conjecture that on the basis of data privatized this way the same rates of convergence as derived below should be attainable.
However, since our paper is motivated by the paper \citep{BeGyWa21} where Laplace perturbation was the method of choice we exclusively stick with this method in this the paper.
\end{remark}

\section{Rate of convergence}
\label{Sec:MainRegression}

For a threshold $t>0$, \citep{BeGyWa21} considered the estimator
\begin{align}
\label{eq:def:mtilde}
	\mtilde_n(x) = \frac{\nutilde_n(A_{h,j})}{\mutilde_n(A_{h,j})}\IND_{\{\mutilde_n(A_{h,j})\ge t\}}\IND_{\{j\le N_n\}} \quad \text{when } x \in A_{h,j},
\end{align}
where
\begin{equation*}
	\nutilde_n(A_{h,j}) = \frac{1}{n} \sum_{i=1}^n Z_{ij} \quad \text{and} \quad \mutilde_n(A_{h,j}) = \frac{1}{n} \sum_{i=1}^n W_{ij}.
\end{equation*}
In \citep{BeGyWa21} the convergence rate $n^{-1/(d+1)}$ was derived (up to a logarithmic term) for Lipschitz continuous functions by specializing \eqref{eq:def:mtilde} with $h=h_n \asymp n^{-1/(2d+2)}$ and $t=t_n \asymp h_n^d/\sqrt{\log n}$. 
However, the proof given in that paper is essentially based on the validity of the strong density assumption (SDA), which means that if $\mu(A_{h,j})>0$, then
\begin{equation}\label{eq:SDA}
    \mu(A_{h,j}) \geq c h^d, \qquad j = 1,\ldots,N_n,
\end{equation}
for some constant $c > 0$.
Moreover, instead of taking \eqref{eq:SDA} as an assumption it was deduced mistakenly from the existence of a density that is lower bounded from below on its support.
Apart from this minor flaw in the proof, imposing \eqref{eq:SDA} is rather artificial since it is not only a condition on the design distribution $\mu$ alone but also on the relationship between the distribution and the sets $A_{h,j}$ of the chosen partition.
Since this condition can hardly be justified in general, it is desirable to eliminate it from the prerequisites.
In general, without the SDA the convergence rate $n^{-1/(d+1)}$ is not attainable via the estimator from \citep{BeGyWa21} (see Remark~\ref{rem:not:optimal} below).

\bigskip

In the sequel, we introduce a novel estimator and bound the private rate of convergence without assuming the SDA.
This disproves the conjecture from \citep{BeGyWa21} where it was guessed that the rate of convergence of any estimate can be arbitrarily slow when the SDA does not hold.

The idea for the general estimator is to include a further modification that in some sense enforces condition \eqref{eq:SDA} to hold (see Remark~\ref{rem:motivation} below for more details).
We again depart from the privatized data \eqref{Eq:Mech1} and \eqref{Eq:Mech2}.
This already guarantees LDP since no further data that depend on the raw data are used in what follows.

In order to define our novel estimator, let $\lambda_n$ denote the uniform distribution on $A_n \defeq \bigcup_{j=1}^{N_n} A_{h,j}$, that is, for any Borel set $A$,
\begin{align*}
\lambda_n(A)=\frac{\lambda(A\cap A_n)}{\lambda(A_n)},
\end{align*}
where $\lambda$ denotes the Lebesgue measure.
We now define our final estimator $[\mhat_n]_T$ by
\begin{align*}
    [\mhat_n]_T(x) \defeq  (-T) \vee (\mhat_n(x) \wedge T)
\end{align*}
where
\begin{align*}
\mhat_n(x) = \frac{\nutilde_n(A_{h,j})}{\muhat_n(A_{h,j})}\IND_{\{\muhat_n(A_{h,j})\ge t\}}\IND_{\{j\le N_n\}} \quad \text{when } x \in A_{h,j},
\end{align*}
and
\begin{align*}
\muhat_n(A_{h,j})&=\frac 34 \mutilde_n(A_{h,j}) + \frac 14 \lambda_n(A_{h,j})\\
&= \left[ \frac 34 \mutilde_n(A_{h,j}) + \frac{1}{4N_n} \right] \IND_{\{j\le N_n\}}.
\end{align*}

The following result states upper risk bounds for the estimators $\mhat_n$ and $[\mhat_n]_T$.
Its proof is deferred to Section~\ref{app:proof:upper}.

\begin{theorem}\label{thm:without:SDA}
Assume that $m$ satisfies \eqref{eq:hoelder} with $0<\beta\leq 1$ and that both $X$ and $Y$ are bounded such that $\lvert X \rvert \leq r$ and $\lvert Y \rvert \leq T$.
Consider the estimator $\mhat_n$ with $t=t_n= \lambda_n(A_{h,1})/8 = 1/(8N_n)$.
Then,
\begin{align}\label{eq:thm:rate}
\EXP \left[ \int ( m(x) - \mhat_n(x) )^{2} \mu(\dd x) \right] 
&\lesssim
\frac{1}{nh_n^d} + h_n^{2\beta} + \frac{\sigma_Z^2}{nt_n^2}\notag\\
&\hspace{1em}+\frac{\sigma_W^2}{nt_n^{2}} + \exp \left( - \frac{8nt_n^2}{9} \right).
\end{align}
As a consequence, taking $h_n \asymp \max\{ (1/(n\alpha^2))^{1/(2\beta +2d)}, (1/n)^{1/(2\beta + d)} \}$ and $\sigma_W$, $\sigma_Z$ as in \eqref{eq:choice:sigmas}, yields
\begin{align}\label{eq:thm:consequence}
\EXP &\left[ \int ( m(x) - [\mhat_n]_T(x) )^{2} \mu(\dd x) \right]\lesssim  (n^{-\frac{2\beta}{2\beta + d}} \vee (n\alpha^2)^{-\frac{\beta}{\beta+d}}) \wedge 1
\end{align}
where the numerical constant hidden in the $\lesssim$ notation depends on $r,T$, and $d$.
\end{theorem}

\begin{remark}\label{rem:not:optimal}
The estimate $\mtilde_n$ defined by (\ref{eq:def:mtilde}) cannot achieve the rate of convergence in \eqref{eq:thm:consequence} without assuming the SDA. 
In order to see this, consider the particular case when there is no privatization and we are given constant, noiseless observations, that is, $Y=C$ a.s.\,for a constant $C \neq 0$.
Then, the estimator $\mtilde_n$ satisfies
\begin{align*}
&\int \EXP [( m(x) -  \mtilde_n(x))^2 ] \mu(\dd x)\\
&=
C^2\sum_{j=1}^{N_n} \EXP [\IND_{ \{ \mu_n(A_{h,j}) < t_n \} }]\mu(A_{h,j})\\
&\ge 
C^2\sum_{j=1}^{N_n} \PROB ( \mu_n(A_{h,j}) < t_n ) \IND_{ \{ \mu(A_{h,j}) < 2t_n \} }\mu(A_{h,j})\\
&\ge 
C^2 \sum_{j=1}^{N_n} \PROB (  \mu_n(A) - \mu(A) < -t_n)  \IND_{ \{ \mu(A_{h,j}) < 2t_n \} }\mu(A_{h,j}).
\end{align*}
Application of Hoeffding's inequality or normal approximation easily shows that the probability in the last line is bounded from below by $1/2$ for $n$ sufficiently large.
Consider now the case $d=1$ and assume that the distribution $\mu$ has a density $f$ on $[0,1]$ satisfying $f(x)=x$ for $0\le x\le \delta$ with some $\delta>0$.
Let the regular partition be given by $[0,h_n),[h_n,2h_n), \ldots$
For $n$ sufficiently large we have $\IND_{ \{ \mu(A_{h,j}) < 2t_n \} } = 1$ at least for those cells contained in $[0,\delta]$ where $xh_n < 2t_n$ holds.

Consequently, with the choice $t_n \asymp h_n/\sqrt{\log n}$ as suggested in \citep{BeGyWa21} we obtain the lower bound
\begin{align*}
\int \EXP [( m(x) -  \mtilde_n(x))^2 ] \mu(\dd x)
&\ge 
\frac{C^2}{2}\int_0^{\lfloor \delta/h_n \rfloor \cdot h_n} \IND_{ \{ xh_n < 2t_n \} }  x\dd x\\
&=
\frac{C^2}{2}\int_0^{\lfloor \delta/h_n \rfloor \cdot h_n} \IND_{ \{ x < 2/\sqrt{\log n} \} }x\dd x\\
&= 
\frac{C^2}{2}\int_0^{2/\sqrt{\log n}} x\dd x\\
&\geq \frac{C^2}{\log n},
\end{align*}
which is a much slower rate than the one in~\eqref{eq:thm:consequence}.
\end{remark}

\begin{remark}\label{rem:motivation}
The rationale behind the construction of the estimator $\mhat_n$ comes from the auxiliary model obtained by replacing the raw data $\mathcal D_n$ with
\begin{equation*}
    \mathcal D\p_n = \{ (X\p_1,Y\p_1),\ldots,(X\p_n,Y\p_n) \}
\end{equation*}
where with probability $3/4$ one has $(X\p_i,Y\p_i)=(X_i,Y_i)$, and with probability $1/4$ one has $X\p_i \sim \lambda_n$ and $Y\p_i = 0$ independently for each $i=1,\ldots,n$.
In this mixture model condition \eqref{eq:SDA} holds by construction.
Recently, in the context of density estimation under local differential privacy the derivation of optimal convergence rates has also been reduced to such a mixture model (see the proof of Proposition~7 in \citep{sart2022density}).
The definition of our estimator $\mhat_n$ is motivated by this approach but our definition does not rely on any additional randomization and replacing $\mutilde$ by $\muhat$ in the definition of $\mhat_n$ may be interpreted as some kind of regularization.
\end{remark}

\section{Lower bound}\label{sec:lower}

In order to prove a lower bound, we restrict ourselves to a specific instance of the general regression model \eqref{def:model}.
This submodel is chosen sufficiently complex to rule out inference with an essentially faster rate than the one obtained in Theorem~\ref{thm:without:SDA}.
More precisely, we consider the regression model with a generic observation $(X,Y) \in \R^d \times \R$ obeying the model
\begin{equation}\label{eq:red:model}
    Y = m(X) + \eta
\end{equation}
where $X$ is distributed uniformly on $[0,1]^d$, the noise $\eta$ is distributed uniformly on $[-1/2,1/2]$, and the regression function $m$ belongs to the Hölder class $\Fclass$ defined as the set of functions satisfying \eqref{eq:hoelder} and with support contained in $[0,1]^d$ and bounded from above by some constant $M >0$.

For the lower bound we allow the potential privacy mechanism to belong to the class $\Qc_\alpha$ of \emph{sequentially interactive} privacy mechanism that generalizes the class of non-interactive mechanism introduced in Section~\ref{s:privacy}.
More precisely, given any ordering of the raw data $X_i$, privatized data $Z_i$ are generated according to conditional probability measures
$Q_i(\boldsymbol{\cdot}\,|(X_i,Y_i)=(x_i,y_i),Z_1=z_1,\ldots,Z_{i-1}=z_{i-1})$, that is, the value $Z_i$ does not only depend on $X_i$ but also on privatized data created by other data holders before.
In this more general case, condition~\eqref{eq:LDP:noninteractive} is replaced by 
\begin{equation*}
    \frac{Q_i(A|(X_i,Y_i)=(x,y),Z_{i-1}=z_{i-1},\ldots,Z_1=z_1)}{Q_i(A|(X_i,Y_i)=(x\p,y\p),Z_{i-1}=z_{i-1},\ldots,Z_1=z_1)} \leq \exp(\alpha)
\end{equation*}
which must be satisfied for any $A \in \Zs_i$, $z_j \in \Zc_j$ for $j=1,\ldots,i-1$,
and all potential values $(x,y), (x\p,y\p)$ of the raw data.

For this general setup, we can prove the following lower bound result. 
Its proof, which is based on Assouad's lemma, is given in Section~\ref{app:proof:lower}.

\begin{theorem}\label{thm:lower}
    It holds
    \begin{align*}
        \inf_{\substack{\mtilde \\ Q \in \Qc_\alpha}} \sup_{\substack{m \in \Fclass}} \EXP \left[ \int (\mtilde(x) - m(x))^2 \dd x \right] \gtrsim (n(e^\alpha - 1)^2)^{-\beta/(\beta + d)} \wedge 1.
    \end{align*}
    where the supremum is taken over all admissible regression functions from the Hölder class $\Fclass$ and the infimum is taken over all estimators $\mtilde$ based on a private sample of size $n$ with raw data from model \eqref{eq:red:model} and all, potentially sequentially interactive, privacy mechanisms $Q \in \Qc_\alpha$. 
\end{theorem}
Combining Theorem~\ref{thm:lower} with the fact that the convergence rate under LDP cannot be faster than then non-private rate $n^{-2\beta/(2\beta +d)}$, it turns out that the rate of convergence derived in Theorem~\ref{thm:without:SDA} is essentially optimal.

\begin{remark}\label{rem:design}
    The lower bound is established for the specific design measure $\mu$ given by the uniform distribution on $[0,1]^d$.
    The proof can be easily adapted to design measures $\mu$ with support $[0,1]^d$ and Lebesgue density bounded away from zero.
\end{remark}
\begin{remark}\label{rem:error}
    Similarly, the theorem is established for the special error distribution given by the uniform distribution on $[-1/2,1/2]$.
    This choice permits explicit calculations of certain total variation distances needed in the proof of the lower bound.
    It is an interesting open question whether the same lower bound for private estimation holds true already when the response variable in the raw data is noise-free, that is, $Y_i = m(X_i)$ for $i=1,\ldots,n$, since this would rule out the possibility of error distributions with a faster convergence rate.
    Note that our proof of Theorem~\ref{thm:lower} does not apply in this case.
\end{remark}

\section{Proofs}
\label{Sec:Proofs}

\subsection{Proof of Theorem~\ref{thm:without:SDA} (Upper bound)}\label{app:proof:upper}

Loosely speaking, the proof of the upper bound decomposes the overall risk into three terms that are bounded separately. The first term (estimate~\eqref{eq:bound:1} below) captures the privatization of the response variable in \eqref{Eq:Mech1} and its upper bound contains the noise level $\sigma_Z$.
The second term (estimate~\eqref{eq:bound:nonprivate} below) yields the classical bound that holds already for the raw data without privatization, and the third and last term (estimate~\eqref{eq:bound:3} below) is the contribution due to the privatization of the covariate values in \eqref{Eq:Mech2} which contains the noise level $\sigma_W$.

\bigskip

We start the proof with the decomposition
\begin{equation*}
	\mhat_n = \mphat_1 + \mphat_2
\end{equation*}
where, for $x \in A_{h,j}$, we set
\begin{equation*}
	\mphat_1(x) = \frac{\frac{\sigma_Z}{n} \sum_{i=1}^n \epsilon_{ij}}{\muhat_n(A_{h,j})} \IND_{\{\muhat_n(A_{h,j})\geq t_n\}} \IND_{\{ j \leq N_n\}}
\end{equation*}
and
\begin{equation*}
	\mphat_2(x) = \frac{ \nu_n(A_{h,j})}{\muhat_n(A_{h,j})} \IND_{\{\muhat_n(A_{h,j})\geq t_n\}} \IND_{\{ j \leq N_n\}}
\end{equation*}
(recall that we assume that $\lvert Y \rvert \leq T$ which implies that $[Y_i]_T = Y_i$ for all $i=1,\ldots,n$).
With this notation, proving \eqref{eq:thm:rate} reduces to show that
\begin{align}\label{eq:bound:1}
\int \EXP[ (\mphat_1(x))^2] \mu(\dd x) 
\leq \frac{\sigma_Z^2}{nt_n^{2}}
\end{align}
and
\begin{align}\label{eq:bound:2}
\int \EXP[ (m(x) - \mphat_2(x))^2 ] \mu(\dd x) 
&\lesssim \frac{1}{nh_n^d} +h_n^{2\beta}\notag\\
&\hspace{1em}+\frac{\sigma_W^2 \vee \sigma_Z^2}{nt_n^{2}} + \exp \left( - \frac{8nt_n^2}{9} \right).
\end{align}

\noindent\emph{Proof of \eqref{eq:bound:1}}: One has
\begin{align*}
\int \EXP[ (\mphat_1(x))^2] \mu(\dd x)  
&= \sum_{j=1}^{N_n} \EXP\left[ \frac{(\frac{\sigma_Z}{n} \sum_{i=1}^n \epsilon_{ij})^2}{(\muhat_n(A_{h,j}))^2} \IND_{\{\muhat_n(A_{h,j})\geq t_n\}} \right] \mu(A_{h,j})\\
&\le 
\frac{\sigma_Z^2}{nt_n^2} \sum_{j=1}^{N_n} \mu(A_{h,j})\\
&\leq
\frac{\sigma_Z^2}{nt_n^2}.
\end{align*}

\noindent \emph{Proof of \eqref{eq:bound:2}:}
Let $\mnpr$ be the modification of $m_n$ where $\mu_n$ is replaced by $\munp= \frac 34 \mu_n + \frac 14 \lambda_n$.
Then, following along the lines of the proof of Theorem~4.3 in \citep{GyKoKrWa02} one can show that
\begin{equation}\label{eq:bound:nonprivate}
\EXP \left[ \int (m(x) - \mnpr (x))^2  \mu(\dd x)  \right] 
\leq \frac{C_1T^2}{nh_n^d} + C_2h_n^{2\beta}
\end{equation}
where $C_1=C_1(d,\mu)$ (more precisely, this constant depends on the measure $\mu$ only through its support) and $C_2=C_2(d)$.
In order to show \eqref{eq:bound:2} it is sufficient to show that
\begin{equation}\label{eq:bound:3}
\EXP \left[ \int (\mphat_2(x) - \mnpr (x))^2  \mu(\dd x)  \right] \lesssim \exp \left( - \frac{8nt_n^2}{9} \right) +  \frac{\sigma_W^2}{nt_n^{2}}.
\end{equation}
In order to prove this bound, note that
\begin{align}\label{eq:bound:4}
\int (\mphat_2(x) - \mnpr (x))^2  \mu(\dd x) \leq J_n
+ 
\frac{16T^2}{9} \mu\left(\R^d \setminus A_n \right),
\end{align}
where
\begin{equation*}
	J_n = \sum_{j=1}^{N_n} (\nu_n(A_{h,j}))^2 \left( \frac{1}{\munp (A_{h,j})} - \frac{1}{\muhat_n(A_{h,j})} \IND_{ \{ \muhat_n(A_{h,j}) \geq t_n \} } \right)^2 \mu(A_{h,j}).
\end{equation*}
Since $\lvert X\rvert \leq r$, the support of $\mu$ is contained in $A_n$, and consequently the second term on the right-hand side of \eqref{eq:bound:4} vanishes.
Therefore it is sufficient to find a bound for $\EXP[J_n]$.
Using that $3\mu_n(A_{h,j})/4 \leq \munp(A_{h,j})$,
\begin{align*}
J_n 
&\leq T^2 \sum_{j=1}^{N_n} (\mu_n(A_{h,j}))^2 \left( \frac{1}{\munp(A_{h,j})} - \frac{1}{\muhat_n(A_{h,j})} \IND_{ \{ \muhat_n(A_{h,j}) \geq t_n \} } \right)^2 \mu(A_{h,j})\\
&\leq 
\frac{16T^2}{9} \sum_{j=1}^{N_n} \left( 1 - \frac{\munp(A_{h,j})}{\muhat_n(A_{h,j})} \IND_{ \{ \muhat_n(A_{h,j}) \geq t_n \} } \right)^2 \mu(A_{h,j})\\
&= \frac{16T^2}{9} \sum_{j=1}^{N_n} \IND_{ \{ \muhat_n(A_{h,j}) < t_n \} }  \mu(A_{h,j})\\
&\qquad
+ \frac{16T^2}{9} \sum_{j=1}^{N_n} \left( 1 - \frac{\munp(A_{h,j})}{\muhat_n(A_{h,j})} \right)^2 \IND_{ \{ \muhat_n(A_{h,j}) \geq t_n \} } \mu(A_{h,j}).
\end{align*} 
Therefore,
\begin{equation}\label{eq:bound:E:J_n}
	\EXP[J_n] \leq \EXP [ J_{n,1} ] + \EXP [ J_{n,2} ]
\end{equation}
where, setting $\mup=\frac{3}{4} \mu + \frac{1}{4} \lambda_n$,
\begin{align*}
J_{n,1} &= \frac{64T^2}{27} \sum_{j=1}^{N_n} \IND_{ \{ \muhat_n(A_{h,j}) < t_n \}}\mup(A_{h,j}), \quad \text{and}\\
J_{n,2} &= \frac{16 T^2}{9} \sum_{j=1}^{N_n} 
\left( 1 - \frac{\munp(A_{h,j})}{\muhat_n(A_{h,j})} \right)^2 \IND_{ \{ \muhat_n(A_{h,j}) \geq t_n \} } \mu(A_{h,j}).
\end{align*}
In order to deal with the expectation of $J_{n,1}$, note that $\mup(A_{h,j})\geq \lambda_n(A_{h,j})/4 = 1/(4N_n) \geq 2t_n$ holds for our choice of $t_n$ by definition.
Thus,
\begin{align*}
\PROB(\muhat_n(A_{h,j}) < t_n) 
&=
\PROB(\muhat_n(A_{h,j}) < t_n,\mup(A_{h,j})\geq 2t_n )\\ 
&\leq 
\PROB( \mup(A_{h,j})-\muhat_n(A_{h,j})\geq t_n )\\
&\leq P_{1,j} + P_{2,j}
\end{align*}
where
\begin{align*}
P_{1,j} &= \PROB\left(  \mu(A_{h,j})-\mu_n(A_{h,j}) \geq \frac{2t_n}{3} \right),\\
P_{2,j} &= \PROB \left( \frac{\sigma_W}{n} \sum_{i=1}^n \zeta_{ij} \geq \frac{2t_n}{3} \right).
\end{align*}
Applying Hoeffding's inequality in the formulation taken from \citep{boucheron2013concentration}, Theorem~2.8, yields
\begin{equation*}
    P_{1,j} \leq \exp \left( - \frac{8nt_n^2}{9} \right),
\end{equation*}
whereas Chebyshev's inequality implies that
\begin{equation*}
    P_{2,j} \leq \frac{9\sigma_W^2}{4nt_n^2}.
\end{equation*}
Hence,
\begin{equation*}
\EXP [J_{n,1}] \leq \frac{64T^2}{27} \exp \left( - \frac{8nt_n^2}{9} \right)  + \frac{16T^2\sigma_W^2}{3nt_n^2} 
\end{equation*}
Furthermore, we have
\begin{align*}
\EXP [ J_{n,2} ]
&=
\frac{16T^2}{9} \sum_{j=1}^{N_n} \EXP\left\{
\left( \frac{\muhat_n(A_{h,j})-\munp(A_{h,j})}{\muhat_n(A_{h,j})} \right)^2 \IND_{ \{ \muhat_n(A_{h,j}) \geq t_n \} }\right\} \mu(A_{h,j})\\
&\leq
\frac{16T^2}{9t_n^2} \sum_{j=1}^{N_n} \EXP [
\left(\muhat_n(A_{h,j})-\munp(A_{h,j})\right)^2 ] \mu(A_{h,j})\\
&\leq 
\frac{16 T^2\sigma_W^2}{9nt_n^2}.
\end{align*}
Putting the bounds obtained for $\EXP [ J_{n,1} ]$ and $\EXP [ J_{n,2} ]$ into \eqref{eq:bound:E:J_n} yields \eqref{eq:bound:3} which proves \eqref{eq:thm:rate}.
\eqref{eq:thm:consequence} follows from \eqref{eq:thm:rate} by taking into account that
\begin{equation*}
	t_n \asymp \lambda_n(A_{h,1}) = h_n^d/\lambda(A_n) \asymp h_n^d/r_n^d
\end{equation*}
and that the term of order $\exp (-8nt_n^2/9)$ in \eqref{eq:thm:rate} is negligible.
The truncation in the definition of $[\mhat_n]_T$ guarantees that the risk is at least bounded by a constant depending on $T$, $r$, and the dimension $d$.

\subsection{Proof of Theorem~\ref{thm:lower} (Lower bound)}\label{app:proof:lower}

The overall strategy to establish the stated private lower bound is similar to the one for the classical lower bound that holds for estimators defined in terms of the raw data. Indeed, the reduction to the pairwise comparison of certain hypotheses parameterized through the corners of a high-dimensional hypercube and the construction of these hypotheses below is borrowed from Chapter~2.6.1 in \citep{tsybakov2009introduction}.
In order to apply Assouad's lemma under privacy constraints one has to use a suitable bound for the Kullback-Leibler divergence of the privatized data under different hypothetical regression functions.
Such a bound will be derived by combining the information theoretical inequality \eqref{eq:bound:Kullback} with a bound for the (squared) total variation distance of the raw data. The use of inequality \eqref{eq:bound:Kullback} has turned out to be essential for proving lower bounds in the density estimation problem before \citep{Butucea2020,Gyoerfi2023}.

\bigskip

We start by introducing some notation that will be used throughout the proof.
Let $K_0 \colon \R \to [0,\infty)$ be a $C^\infty$-function such that \eqref{eq:hoelder} is satisfied with constant equal to $1$, $\lVert K_0\rVert_\infty \leq 1$, and $\supp(K_0) \subseteq [0,1]$.
For $x = (x_1,\ldots,x_d) \in [0,1]^d$, define the function $K\colon [0,1]^d\to \R$ via $K(x) = \min_{i=1,\ldots,d} K_0(x_i)$.
We restrict the complexity of the whole problem by restricting ourselves to a finite set of hypotheses parameterized by $\theta \in \Theta \defeq \{ 0,1 \}^{k^d}$ for some positive integer $k$ that will be specified below.
Set $c = C \wedge M$.
For any $j = (j_1,\ldots,j_d) \in \{ 0,\ldots,k-1 \}^d$ define the function $K_j$ by
\begin{equation*}
    K_j(x) = ck^{-\beta} K(kx_1-j_1,\ldots,kx_d-j_d).
\end{equation*}
It is readily checked that $\supp(K_j) \subseteq B_j \defeq \times_{i=1}^d [j_i/k,(j_i+1)/k]$.
For any $\theta \in \Theta$, we consider the candidate regression function
\begin{equation*}
    m_\theta = \sum_{j} \theta_j K_j,
\end{equation*}
where the sum is taken over all multi-indices $j \in \{ 0,\ldots,k-1 \}^d$.
By construction $m_\theta$ belongs to $\Fclass$ for any $\theta \in \Theta$.
Let us now assume that the raw data have been privatized by means of an arbitrary privacy mechanism $Q \in \Qc_\alpha$, and let $\mtilde$ be any estimator defined in terms of the outcome $Z$ of $Q$.
We denote with $\PROB_\theta$ the distribution of the tupel $(X_1,Y_1)$ and with $Q\PROB_\theta^n$ the distribution of $Z=(Z_1,\ldots,Z_n)$ when the true regression function is $m_\theta$.
We also write $\EXP_\theta$ for the expectation operator in this case.

After these preliminaries we start the proof with the observation that for any $\theta \in \Theta$,
\begin{align*}
    \EXP_\theta \left[ \int_{[0,1]^d} (\mtilde(x) - m_\theta(x))^2 \dd x \right] &= \sum_j \EXP_\theta \left[ \int_{B_j} (\mtilde(x) - m_\theta(x))^2 \dd x \right]\\
    &= \sum_j \EXP_\theta [\rho_j^2(\mtilde,\theta_j)],
\end{align*}
where
\begin{equation*}
    \rho_j(\mtilde,\theta_j) = \left( \int_{B_j} (\mtilde(x) - \theta_j K_j(x))^2 \dd x \right)^{1/2}.
\end{equation*}
Putting $\thetahat_j = \argmin_{t \in \{ 0,1 \}} \rho_j(\mtilde,t)$, we have
\begin{equation*}
    \rho_j(\mtilde, \theta_j) \geq \frac{\lVert K_j \rVert_{2}}{2} \cdot \lvert \thetahat_j - \theta_j \rvert.
\end{equation*}
Hence, using that $\lVert K_j \rVert_2^2 = c^2k^{-2\beta-d} \lVert K \rVert_2^2$, we obtain
\begin{align*}
    \EXP_\theta \left[ \int_{[0,1]^d} (\mtilde(x) - m_\theta(x))^2 \dd x \right] &\geq \frac {c^2 \lVert K \rVert_2^2}{4} k^{-2\beta-d} \EXP_\theta[\rho(\thetahat,\theta)]
\end{align*}
where $\thetahat = (\thetahat_j)$ and $\rho(\theta,\thetap)= \sum_j \IND_{\{\theta_j \neq \theta_j^\prime\}}$ denotes the Hamming distance between $\theta$ and $\thetap$.
Consequently,
\begin{align*}
    \sup_{\substack{m \in \Fclass}} &\EXP \left[ \int_{[0,1]^d} (\mtilde(x) - m(x))^2 \dd x \right]\\
    &\geq \sup_{\theta \in \Theta} \EXP_\theta \left[ \int_{[0,1]^d} (\mtilde(x) - m_\theta(x))^2 \dd x \right]\\
    &\geq  \frac{c^2 \lVert K \rVert_2^2}{4} k^{-2\beta-d} \inf_{\thetahat} \sup_{\theta \in \Theta} \EXP_\theta [\rho(\thetahat,\theta)].
\end{align*}
In order to bound the quantity $\inf_{\thetahat} \sup_{\theta \in \Theta } \EXP_\theta [\rho(\thetahat,\theta)]$, we use Statement~(iv) of Theorem~2.12 in \citep{tsybakov2009introduction} which relies on a finite bound on the Kullback-Leibler distance $K(Q\PROB_\theta^n, Q\PROB_\thetap^n)$ for $\theta, \thetap$ such that $\rho(\theta,\thetap) = 1$.
In order to obtain such a bound, first note that Equation~(14) in \citep{duchi2018} yields
\begin{equation}\label{eq:bound:Kullback}
    K(Q\PROB_\theta^n, Q\PROB_\thetap^n) \leq 4n(e^\alpha - 1)^2 V^2(\PROB_\theta,\PROB_\thetap)
\end{equation}
where $V(\PROB,\mathbf Q)$ denotes the total variation distance between two probability measures.
Thus, it remains to find a bound for $V(\PROB_\theta,\PROB_\thetap)$.
In order to bound this quantity, note that under model \eqref{eq:red:model} the vector $(X,Y)$ has a Lebesgue density $\varphi_m$ that is equal to $1$ on the set
$$\{ (x,y) \in \R^{d+1} : x \in [0,1]^d, y \in [m(x)-1/2,m(x)+1/2] \},$$
and equal to $0$ otherwise.
Now, let $\theta, \thetap$ be such that $\rho(\theta,\thetap) = 1$.
Then, a direct calculation using Scheff{\'e}'s Theorem yields
\begin{align*}
    V(\PROB_\theta,\PROB_\thetap) &= \frac{1}{2} \int_{\R^{d+1}} \lvert \varphi_{m_\theta}(x,y) - \varphi_{m_\thetap}(x,y) \rvert \dd x \dd y\\
    &\leq k^{-d} \lVert K_j \rVert_\infty\\
    &= c k^{-d-\beta} \lVert K \rVert_\infty. 
\end{align*}
Combining this bound with \eqref{eq:bound:Kullback} for $k \asymp (n(e^\alpha - 1))^{1/(2\beta+ 2d)} \vee 1$ yields
\begin{equation*}
    K(Q\PROB_\theta^n, Q\PROB_\thetap^n) \lesssim 1.
\end{equation*}
By application of \citep{tsybakov2009introduction}, Theorem~2.12, Statement~(iv), we obtain
\begin{equation*}
    \sup_{\substack{m \in \Fclass}} \EXP \left[ \int_{[0,1]^d} (\mtilde(x) - m(x))^2 \dd x \right] \gtrsim (n(e^\alpha - 1)^2)^{-\beta/(\beta + d)} \wedge 1,
\end{equation*}
which implies the claim since $Q \in \Qc_\alpha$ and $\mtilde$ were arbitrary.

\section*{Funding information}

The research of Martin Kroll was supported by the German Research Foundation (DFG) under the grant DFG DE 502/27-1. The research of L\'aszl\'o Gy\"orfi has been supported by the National Research, Development and Innovation Fund of Hungary under the 2019-1.1.1-PIACI-KFI-2019-00018 funding scheme.

\printbibliography

\end{document}